\documentclass[reqno,11pt]{amsart}

\usepackage{amssymb}
\usepackage[all,cmtip]{xy}
\usepackage{leftidx}
\newtheorem{MainTheorem}{Theorem}
\newtheorem{Proposition}{Proposition}[section]

\newtheorem{Definition}[Proposition]{Definition}
\newtheorem{Lemma}[Proposition]{Lemma}
\newtheorem{Theorem}[Proposition]{Theorem}

\DeclareMathOperator{\Val}{Val}
\DeclareMathOperator{\vol}{vol}

\DeclareMathOperator{\SO}{SO}

\DeclareMathOperator{\pd}{pd}
\DeclareMathOperator{\nc}{nc}
\newcommand{\ed}{\mathrm{d}}
\newcommand{\D}{\mathrm{D}}
\newcommand{\Q}{\mathrm{Q}}
\newcommand{\Ro}{\mathrm{R}}
\newcommand{\So}{\mathrm{S}}
\DeclareMathOperator{\Curv}{Curv}
\DeclareMathOperator{\Area}{Area}
\DeclareMathOperator{\glob}{glob}
\DeclareMathOperator{\Sym}{Sym}
\DeclareMathOperator{\PD}{PD}
\DeclareMathOperator{\contr}{contr}
\newcommand{\R}{\mathbb{R}}
\newcommand{\C}{\mathbb{C}}
\newcommand\flag[2]{\left[\begin{array}{c} #1\\ #2 \end{array}
  \right]} 
  \newcommand{\largewedge}{\mbox{\Large $\wedge$}}

\title{Dual area measures and local additive kinematic formulas}

\author{Andreas Bernig} 
\address{Institut f\"ur Mathematik, Goethe-Universit\"at Frankfurt am Main, 
Robert-Mayer-Str. 10, 60054 Frankfurt, Germany}
\email{bernig@math.uni-frankfurt.de}

\thanks{Supported by DFG grant BE 2484/5-2.}
\begin{document}


\subjclass[2000]{53C65, 
52A22} 

\keywords{Area measure, valuation, kinematic formula, integral geometry}

\begin{abstract}
We prove that higher moment maps on area measures of a euclidean vector space are injective, while the kernel of the centroid map equals the image of the first variation map.

Based on this, we introduce the space of smooth dual area measures on a finite-dimensional euclidean vector space and prove that it admits a natural convolution product which encodes the local additive kinematic formulas for groups acting transitively on the unit sphere. 
 
As an application of this new integral-geometric structure, we obtain the local additive kinematic formulas in hermitian vector spaces in a very explicit way.
\end{abstract}

\maketitle

\section{Introduction}

\subsection{General background}

Let $V$ be a euclidean vector space of dimension $n$ and $\mathrm{SO}(n)$ its rotation group. We denote by $\overline{\SO(n)}:=\mathrm{SO}(n) \ltimes V$ the euclidean group, endowed with the product of the Haar probability measure and the Lebesgue measure. Then the following kinematic formulas play an important role in classical integral geometry.
\begin{align*}
 \int_{\overline{\SO(n)}} \mu_i(K \cap \bar g L) \ed\bar g & =\flag{n+i}{i} \sum_{k+l=n+i} \flag{n+i}{k}^{-1} \mu_k(K)\mu_l(L),\\
 \int_{\SO(n)} \mu_i(K + gL) \ed g & =\flag{2n-i}{n-i} \sum_{k+l=i}
\flag{2n-i}{n-k}^{-1} \mu_k(K)\mu_l(L).
\end{align*}

In both formulas, $\flag{n}{k}:=\binom{n}{k} \frac{\omega_n}{\omega_k\omega_{n-k}}$ denotes the flag coefficient ($\omega_k$ is the volume of the $k$-dimensional unit sphere) and $\mu_i$ denotes the $i$-th intrinsic volume, which can be defined in a number of equivalent ways. For our purpose we will use the following characterization, due to Hadwiger, of $\mu_i$: it is the unique convex valuation (i.e. $\mu(K \cup L)=\mu(K)+\mu(L)-\mu(K \cap L)$ whenever $K,L,K \cup L$ are compact convex sets) which is continuous with respect to Hausdorff topology, invariant under translations and rotations, $i$-homogeneous (i.e. $\mu_i(tK)=t^i\mu_i(K)$) and equals the $i$-dimensional Lebesgue measure on $i$-dimensional compact convex sets.

The first formula is called {\it intersectional kinematic formula} (or sometimes just {\it kinematic formula}), while the second is called {\it additive kinematic formula} or {\it rotational mean value  formula}. 

It was noted by Nijenhuis \cite{nijenhuis74} that after some rescaling of the $\mu_i$'s, all constants in the intersectional kinematic formula equal $0$ or $1$. He speculated that there may be some algebraic structure behind this observation. Indeed, there is an easy explanation of this fact which motivates our work. For this, let $\Val$ be the vector space of continuous and translation invariant valuations.  McMullen \cite{mcmullen77} showed that 
\begin{displaymath}
\Val=\bigoplus_{\substack{k=0,\ldots,n\\ \epsilon= \pm}} \Val_k^\epsilon,
\end{displaymath}
where $\Val_k^\epsilon=\{\phi \in \Val: \mu(tK)=t^k\mu(K), t >0, \mu(-K)=\epsilon \mu(K)\}$ is the space of $k$-homogeneous valuations which are even/odd. 

By Hadwiger's theorem \cite{hadwiger_vorlesung, klain_rota}, the subspace $\Val^{\SO(n)}$ of rotation invariant elements is spanned by $\mu_0,\ldots,\mu_n$, in particular it is of finite dimension $n+1$. Let us rewrite the kinematic formula in the form of an operator
\begin{align*}
 k:\Val^{\SO(n)} & \to \Val^{\SO(n)} \otimes \Val^{\SO(n)}\\
 \phi & \mapsto \left[(K,L) \mapsto \int_{\overline{\SO(n)}} \phi(K \cap \bar gL)\ed\bar g \right].
\end{align*}
Explicitly, 
\begin{displaymath}
 k(\mu_i)=\flag{n+i}{i} \sum_{k+l=n+i} \flag{n+i}{k}^{-1} \mu_k \otimes \mu_l.
\end{displaymath}
By a change of variables and Fubini's theorem, one easily checks that $k$ is a cocommutative, coassociative coproduct on $\Val^{\SO(n)}$. Hence the adjoint map 
\begin{displaymath}
 k^*: \Val^{\SO(n)*} \otimes \Val^{\SO(n)*} \to \Val^{\SO(n)*} 
\end{displaymath}
makes $\Val^{\SO(n)*}$ into a commutative, associative algebra with respect to the adjoint map $k^*$. By looking at degrees of homogeneity, one sees that this algebra is isomorphic to the polynomial algebra $\R[t]/(t^{n+1})$. In the basis given by $1,t,\ldots,t^n$, each structure constant is obviously equal to $0$ or $1$, which explains Nijenhuis' observation. A similar reasoning also applies to the additive kinematic formulas 
\begin{align*}
 a:\Val^{\SO(n)} & \to \Val^{\SO(n)} \otimes \Val^{\SO(n)}\\
 \phi & \mapsto \left[(K,L) \mapsto \int_{\SO(n)} \phi(K +gL)\ed g \right],
\end{align*}
which turn $\Val^{\SO(n)*}$ into an algebra, which is also isomorphic to $\R[t]/(t^{n+1})$.

Both types of kinematic formulas admit generalizations to other groups as follows. For this, let $G$ be a closed subgroup of $\SO(n)$ and denote by $\Val^G$ the vector space of continuous, translation and $G$-invariant valuations. Alesker has shown that $\dim \Val^G$ is finite-dimensional if and only if $G$ acts transitively on the unit sphere. The compact connected groups which act transitively and effectively on some sphere are known to belong to one of the series
\begin{displaymath}
 \mathrm{SO}(n), \mathrm{U}(n), \mathrm{SU}(n), \mathrm{Sp}(n), \mathrm{Sp}(n) \cdot \mathrm{U}(1),
\mathrm{Sp}(n)\cdot \mathrm{Sp}(1),
\end{displaymath}
or equal one of the three exceptional groups 
\begin{displaymath}
 \mathrm{G}_2, \mathrm{Spin}(7), \mathrm{Spin}(9). 
\end{displaymath}
For simplicity, we will call a group from this list a {\it transitive group}.

In each of these cases $\Val^{G*}$ becomes an algebra with respect to $k_G^*$ (the adjoint of the kinematic operator $k_G$ which is defined by replacing $\overline{\SO(n)}$ by $\overline{G}=G \ltimes V$), and also with respect to $a_G^*$. A fundamental fact of integral geometry is that all these algebras are subalgebras of the space of smooth translation invariant valuations endowed with the Alesker product \cite{alesker04_product} in the case of intersectional kinematic formulas and endowed with the convolution product \cite{bernig_fu06} in the case of additive kinematic formulas. 
More precisely, $(\Val^{sm},\cdot)$ satisfies a version of Poincar\'e duality, so that there is a map $\PD:\Val^{sm} \to \Val^{sm,*}$. Similarly, $(\Val^{sm},*)$ satisfies a version of Poincar\'e duality, giving rise to another map $\PD: \Val^{sm} \to \Val^{sm,*}$ (since these two maps agree up to a sign, which equals $1$ in the cases to be considered here, we will not distinguish by our notation). Then the following diagram commutes
\begin{displaymath}
 \xymatrix{\Val^{sm} \otimes \Val^{sm} \ar[r] \ar[d]^{\PD \otimes \PD} & \Val^{sm} \ar[d]^{\PD}\\ \Val^{sm,*} \otimes \Val^{sm,*} \ar[r] \ar[d] & \Val^{sm,*}\ar[d]  \\\Val^{G*} \otimes \Val^{G*} \ar[r] & \Val^{G*}
 }
\end{displaymath}
where the vertical maps in the lower square are dual to the inclusion maps $\Val^G \to \Val^{sm}$, and the upper horizontal map is the Alesker product in the case of intersectional kinematic formulas and the convolution product in the case of additive kinematic formulas. As a side remark, we note that the intersectional and the additive kinematic formulas are related by the Alesker-Fourier transform, see \cite{bernig_fu06}.

In the case $G=\mathrm U(n)$, the kinematic formulas as well as the additive kinematic formulas were obtained in \cite{bernig_fu_hig}. The algebra structure was computed earlier by Fu \cite{fu06}, who showed that $\Val^{\mathrm U(n)}$ is isomorphic to $\R[t,s]/\langle f_{n+1},f_{n+2}\rangle$, where $\log(1+tx+sx^2)=\sum_{k=0}^\infty f_k(t,s)x^k$. For the integral geometry of the other transitive groups, we refer to \cite{bernig_aig10} and the references therein.

The kinematic formulas admit localized versions, which apply to {\it smooth curvature measures} in the case of intersectional kinematic formulas and to {\it smooth area measures} in the case of additive kinematic formulas. A smooth curvature measure is a valuation with values in the space of signed measures on $V$, while a smooth area measure is a valuation with values in the space of signed measures on the unit sphere. The technical notion of smoothness which will be recalled in Section \ref{sec_smooth_valuations}. 

We let $\Curv$ (resp. $\Area$) denote the spaces of smooth, translation covariant curvature measures (area measures resp.). If $G$ is a transitive group, then $\dim \Curv^G<\infty, \dim \Area^G<\infty$. If $\Phi_1,\ldots,\Phi_N$ is a basis of $\Curv^G$, then the local kinematic formulas are given by 
\begin{displaymath}
 \int_{\bar G} \Phi_i(K \cap \bar gL,\kappa \cap \bar g \lambda) \ed\bar g=\sum_{k,l} \tilde c_{k,l}^i \Phi_k(K,\kappa) \Phi_l(L,\lambda),
\end{displaymath}
where $\kappa,\lambda$ are bounded Borel subsets of $V$. The existence of such formulas was shown by Fu \cite{fu90}, and as above we obtain a cocommutative coassociative coproduct $K_G$ on $\Curv^G$, or equivalently an algebra structure $K_G^*$ on $\Curv^{G*}$. In the case $G=\mathrm U(n)$, this algebra structure was recently obtained as follows, see \cite{bernig_fu_solanes, bernig_fu_solanes_proceedings}.

Define polynomials $f,p,q$ by 
\begin{align*}
 \log(1+tx+sx^2)& =\sum_{k=0}^\infty f_k(t,s)x^k\\
 \frac{1}{1+tx+sx^2} & = \sum_{k=0}^\infty p_k(t,s)x^k\\
 -\frac{1}{(1+tx+sx^2)^2} & = \sum_{k=0}^\infty q_k(t,s)x^k.
\end{align*}

\begin{Theorem}[\cite{bernig_fu_solanes_proceedings}]
There is an algebra isomorphism
\begin{displaymath}
\Curv^{U(n)*} \cong \C[t,s,v]/\langle f_{n+1},f_{n+2},q_{n-1}v,q_nv, (v+t(4s-t^2))^2 \rangle. 
\end{displaymath}
\end{Theorem}

Similarly, as shown by Wannerer \cite[Theorem 2.1]{wannerer_unitary_module}, if $\Psi_1,\ldots,\Psi_m$ is a basis of $\Area^G$, then there are local additive kinematic formulas
\begin{displaymath}
 \int_{G} \Psi_i(K +gL,\kappa \cap g \lambda) \ed g=\sum_{k,l} \hat c_{k,l}^i \Psi_k(K,\kappa) \Psi_l(L,\lambda),
\end{displaymath}
where $\kappa,\lambda$ are Borel subsets of the unit sphere. We obtain again some operator $A_G$, whose adjoint turns $\Area^{G*}$ into an algebra which was computed by Wannerer \cite{wannerer_area_measures, wannerer_unitary_module} in the case $G=\mathrm U(n)$. 

\begin{Theorem}[Wannerer]
 There is an algebra isomorphism  
 \begin{displaymath}
  \Area^{\mathrm U(n)*} \cong \R[t,s,v]/\langle f_{n+1}(t,s),f_{n+2}(t,s),p_n(t,s)-q_{n-1}(t,s)v,v^2\rangle.
 \end{displaymath}
\end{Theorem}

The main ideas of the proof are relevant for the present paper, so we will describe them briefly here. First, there are two important maps, called {\it first variation map} and {\it moment map}. The first variation map is a map 
\begin{displaymath}
 \delta:\Val^{sm} \to \Area
\end{displaymath}
which is uniquely defined by the property 
\begin{displaymath}
 \left.\frac{d}{dt}\right|_{t=0} \mu(K+tL)=\int_{S^{n-1}} h_L \ed(\delta \phi(K)),
\end{displaymath}
where $h_L$ is the support function of $L$. We will write $\delta_k:\Val^{sm}_{k+1} \to \Area_k, 0 \leq k \leq n-1$ for the restriction of $\delta$ to $\Val^{sm}_{k+1}$. The $r$-th moment map is the map
\begin{displaymath}
 M^r:\Area \to \Val \otimes \Sym^rV, \quad \Phi \mapsto \left[K \mapsto \int_{S^{n-1}} \vec{y}^r \ed(\Phi(K))(\vec y)\right].
\end{displaymath}
Here $\vec{y}^r \in \Sym^r V$ stands for the $r$-th symmetric power of the vector $\vec{y} \in S^{n-1}$.

The zero-th moment map is called {\it globalization map} and denoted by 
\begin{displaymath}
\glob:\Area \to \Val. 
\end{displaymath}
The first moment map $M^1$ is called {\it centroid map}. It is easy to see that the image of $\delta$ belongs to the kernel of $M^1$. In the unitary case, Wannerer has shown that the kernel of $M^1$ equals the image of $\delta$. 

Note that, when applying the globalization map to the local additive kinematic formulas on both sides, we obtain the global additive kinematic formulas. More generally, the image of $M^r$ is a tensor-valued valuation, and for such valuations there exist global additive kinematic formulas \cite{bernig_hug}. Wannerer has shown that the additive kinematic formulas for area measures and for tensor valuations are compatible with respect to the moment maps. It turns out that the second moment map $M^2$ is injective on unitarily invariant area measures, which enables Wannerer to obtain enough information from the tensor case to prove his theorem. 

Although Wannerer's theorem is enough to write down the local additive kinematic formulas, the result is not very explicit; and to work out examples, even in low dimensions, requires some extra work. The full array of local additive kinematic formulas in dimension $2$, and some coefficients in dimension $3$ were found in \cite{wannerer_area_measures}. A more explicit formula, which however only works for the classical surface area measure, was given by Solanes \cite[Theorem 13]{solanes}. Our Theorem \ref{mainthm_kinform_hermitian} may be seen as a generalization of Solanes' formula to all unitarily invariant area measures (even if in the case of the classical surface area our formula is formally different from Solanes' formula). 

To state Solanes' result, let us introduce some notation. In \cite{bernig_fu_hig}, each intrinsic volume $\mu_k$ was written as $\mu_k=\sum_q \mu_{k,q}$, where $\max\{0,k-n\} \leq q \leq \frac{k}{2} \leq n$, with the $\mu_{k,q}$ forming a basis of $\Val^{\mathrm U(n)}$. These elements are called {\it hermitian intrinsic volumes}. 

A basis of $\Area^{\mathrm U(n)}$ was obtained by Park \cite{park02}, see also \cite{bernig_fu_hig}. It consists of elements 
\begin{align*}
 B_{k,q}, & \quad \max\{0,k-n\} \leq q < \frac{k}{2} <n,\\
 \Gamma_{k,q}, & \quad \max\{0,k-n+1\} \leq q \leq \frac{k}{2}<n.
\end{align*}
which are defined in terms of differential forms on the sphere bundle of $\C^n$, see Section \ref{sec_hermitian_case}. The classical surface area measure is in this notation $S_{2n-1}:=B_{2n-1,n-1}$. 
They satisfy $\glob B_{k,q}=\glob \Gamma_{k,q}=\mu_{k,q}$. In particular, the kernel of $\glob$ is spanned by the area measures 
\begin{displaymath}
 N_{k,q}:=\frac{2(n-k+q)}{2n-k}(\Gamma_{k,q}-B_{k,q}), \quad \max\{0,k-n+1\} \leq q < \frac{k}{2} <n.
\end{displaymath}
We let $B$ denote the span of the $B_{k,q}$'s and $\Gamma$ the span of the $\Gamma_{k,q}$'s.

Solanes used a slightly different notation and denoted $\Gamma_{2q,q}$ by $B_{2q,q}$. The $B_{k,q}$'s and $N_{k,q}$'s form a basis of $\Area^{\mathrm U(n)}$ and the first variation map may be decomposed as $\delta=\delta_B+\delta_N$, where $\delta_B,\delta_N$ are the projections onto the corresponding spaces. These maps can be written down in a very explicit way but we will not need it here. Solanes also defined the map $l_B:\Val^{\mathrm U(n)} \to \Area^{\mathrm U(n)}$ such that $l_B(\mu_{k,q})=B_{k,q}$. 

\begin{Theorem}[Solanes]
\begin{displaymath}
 A_{\mathrm U(n)}(S_{2n-1})=(\delta \otimes l_B+l_B \otimes \delta_N) a_{\mathrm U(n)}(\vol).
\end{displaymath}
\end{Theorem}

This shows that the local additive kinematic formula for the surface area measure can be obtained from the global kinematic formula by a very explicit and easy operation. 

\subsection{Results of the present paper}

Our first main theorem is of independent interest and will be a key element in the proof of the other main theorems. It generalizes \cite[Theorems 2.19 and 4.8]{wannerer_unitary_module}, where it is proved under the additional assumption of $\mathrm U(n)$-invariance. 

\begin{MainTheorem} \label{mainthm_moments}
\begin{enumerate}
 \item The kernel of the first moment map $M^1:\Area_k \to \Val_k \otimes V$ equals the image of the first variation map $\delta_k:\Val_{k+1} \to \Area_k$ for each $0 \leq k \leq n-1$.
 \item The higher moment maps $M^r:\Area_k \to \Val_k \otimes \Sym^rV, r \geq 2, 1 \leq k \leq n-1$ are injective. 
\end{enumerate}
\end{MainTheorem}

Note that $\Area_0$ is simply the space of smooth signed measures and $M^r$, the $r$-th moment map, is not injective on this space, which is why we have to assume that $k \geq 1$ in Part (2).

The proof of this theorem is the technical part of the paper and uses a careful analysis of differential forms.

Let us now describe an application of this theorem in integral geometry. 

Let $SV=V \times S^{n-1}$ denote the sphere bundle of $V$ and $\pi=\pi_1:SV \to V$ the projection to the first factor, $\pi_2:SV \to S^{n-1}$ the projection to the second factor. By $\Omega^{k,l}:=\Omega^{k,l}(SV)^{tr}$ we will denote the space of translation invariant forms of bidegree $(k,l)$ on $SV$ and by $\Omega^k$ the space of translation invariant forms of total degree $k$. Given $\rho \in \Omega^{2n-1}$ we write $\int \rho$ for the real number $c$ such that $\pi_*\rho=c \vol$, where $\vol$ is the volume form on $V$.

\begin{Definition}
 Let the space $\Area$ of smooth area measures be endowed with the quotient topology of the usual Fr\'echet topology on the space of $(n-1)$-forms on $SV$ (see Section \ref{sec_smooth_valuations} for details). Elements of the dual space $\Area^*$ will be called {\it dual area measures}. A dual area measure $L$ is called smooth if there exists some $n$-form $\tau$ such that 
 \begin{displaymath}
  \langle L,\Phi\rangle=\int \omega \wedge \tau,
 \end{displaymath}
 whenever $\omega \in \Omega^{n-1}$ and $\Phi$ is the area measure induced by $\omega$. The space of smooth dual area measures is denoted by $\Area^{*,sm}$.
\end{Definition}

Let $G$ be a transitive group. The inclusion $\Area^G \subset \Area$ induces a projection $q_G:\Area^{*,sm} \to \Area^{G*}$. Let $A:\Area^G \to \Area^G \otimes \Area^G$ be the local additive kinematic formula and $A^*:\Area^{G*} \otimes \Area^{G*} \to \Area^{G*}$ its adjoint.

\begin{MainTheorem} \label{main_thm_ftaig}
There exists a natural convolution product $*$ on the space $\Area^{*,sm}$ of smooth dual area measures such that for each transitive group $G$ the following diagram commutes
\begin{displaymath}
 \xymatrix{
 \Area^{*,sm} \otimes \Area^{*,sm} \ar[d]^{q_G \otimes q_G} \ar[r]^-{*} & \Area^{*,sm} \ar[d]^{q_G}\\
 \Area^{G*} \otimes \Area^{G*} \ar[r]^-{A^*} & \Area^{G*}
 }
\end{displaymath}
\end{MainTheorem}

We note that we could also, in the spirit of \cite{bernig_faifman}, introduce a partial convolution product on the whole space $\Area^*$. Since we don't have any application of this more general construction, we will not present the technical details here.

The proof of Theorem \ref{main_thm_ftaig} follows the strategy of \cite{wannerer_area_measures}: we use the moment maps to pass from local additive formulas to additive kinematic formulas for tensor valuations. The fact that higher moment maps are injective will ensure us that there is no loss of information in this process. 

Let us also note that the computation of the convolution product on $\Area^{*,sm}$ is very simple and only involves some algebraic operations (Hodge-star operator and wedge product). 

Let us now consider the hermitian case $G=\mathrm{U}(n)$. As explained above, the global additive kinematic formulas are now well-known. Our next theorem shows that the local additive kinematic formulas follow from the global ones by some simple algebraic operations.

Define operators $l_\Gamma,l_B:\Val^{\mathrm U(n)} \to \Area^{\mathrm U(n)}, J: \Area^{\mathrm U(n)} \to \Area^{\mathrm U(n)}$ by 
\begin{align*}
 l_\Gamma(\mu_{k,q}) & := \begin{cases} \Gamma_{k,q} & \max\{0,k-n+1\} \leq q \leq  \frac{k}{2} \\
 0 & q=k-n,\end{cases}\\
 l_B(\mu_{k,q}) & :=  \begin{cases} B_{k,q} & \max\{0,k-n\} \leq q < \frac{k}{2} \\
 0 & q=\frac{k}{2},\end{cases}\\
  JB_{k,q} & =\frac{2c_{n,k,q}}{c_{n,k-1,q}} \Gamma_{k-1,q}, \quad \max\{0,k-n\} \leq q < \frac{k}{2}<n,\\
  J\Gamma_{k,q} & =- \frac{c_{n,k,q}}{2c_{n,k+1,q}} B_{k+1,q}, \quad \max\{0,k-n+1\} \leq q \leq \frac{k}{2}<n. 
\end{align*}
The constants $c_{n,k,q}$ are defined in \eqref{eq_defcnkq}. Note that $\Area^{\mathrm U(n)}=B \oplus \Gamma$. It is therefore enough to determine the local additive kinematic formulas on $B$ and on $\Gamma$. This is achieved by our next main theorem.

\begin{MainTheorem} \label{mainthm_kinform_hermitian}
The local additive kinematic formulas follow from the global additive kinematic formulas by the following commuting diagram
\begin{displaymath}
 \xymatrix{
  B \ar[r]^-A \ar[d]^J & B \otimes \Gamma \oplus \Gamma \otimes B\\
 \Gamma \ar[d]^\glob \ar[r]^-A & \Gamma \otimes \Gamma \ar[u]_{-(J \otimes \mathrm{id}+\mathrm{id} \otimes J)}\\
 \Val^{\mathrm U(n)} \ar[r]^-a  & \Val^{\mathrm U(n)} \otimes \Val^{\mathrm U(n)} \ar[u]_{l_\Gamma \otimes l_\Gamma}
 }
\end{displaymath}
\end{MainTheorem}

Let us write this down more explicitly. Write the global additive kinematic formulas as 
\begin{displaymath}
 a(\mu_{m,r})=\sum_{k,q,l,p} c_{k,q,l,p}^{m,r} \mu_{k,q} \otimes \mu_{l,p}.
\end{displaymath}
Then for $\max\{0,m-n+1\} \leq r \leq \frac{m}{2}<n$ we have 
\begin{displaymath}
 A(\Gamma_{m,r}) =\sum_{k,q,l,p} c_{k,q,l,p}^{m,r} \Gamma_{k,q} \otimes \Gamma_{l,p}, 
\end{displaymath}
and for $\max\{0,m-n\} \leq r < \frac{m}{2}<n$ we have
\begin{multline*}
A(B_{m,r}) = \frac{c_{n,m,r}}{c_{n,m-1,r}} \sum_{k,q,l,p} c^{m-1,r}_{k,q,l,p} \cdot \\ 
\cdot \left(\frac{c_{n,k,q}}{c_{n,k+1,q}} B_{k+1,q} \otimes \Gamma_{l,p}+\frac{c_{n,l,p}}{c_{n,l+1,p}} \Gamma_{k,q} \otimes B_{l+1,p}\right). 
\end{multline*}

In other words, the formulas for the $\Gamma$ follow by just replacing every occurrence of $\mu$ by $\Gamma$, while the formulas for the $B$'s follow by a slightly more complicated, but elementary operation from the global formulas. The proof of this theorem is based on some elementary consequences of Theorem \ref{main_thm_ftaig}.

\subsubsection{Plan of the paper}

In Section \ref{sec_smooth_valuations} we collect some definitions like smooth valuations and area measures, and some important maps between them like the globalization map, the first variation map and the moment maps.

In Section \ref{sec_rumin_de_rham} we adapt the results from \cite{rumin94} to our situation and prove a formula for the Rumin operator which turns out to be useful.

The technical heart of the paper is contained in Sections \ref{sec_centroid_map} and \ref{sec_higher_moments}, where we study the kernel of the centroid map (i.e. the first moment map) and show that the higher moment maps are injective. For this, we need a careful study of the Rumin operator on tensor-valued forms.  

In Section \ref{sec_dual_area_measures} we introduce the concept of a smooth dual area measure, provide this space with a natural convolution product and show how this new algebraic structure  encodes local additive kinematic formulas. The final section is devoted to the  important case of hermitian integral geometry. We use the new convolution product on smooth dual area measure to obtain the local additive kinematic formulas in this case in a very explicit way.  

\subsubsection*{Acknowledgments}
I thank Gil Solanes for many useful comments on a first version of this paper.

\section{Smooth valuations and area measures}
\label{sec_smooth_valuations}

In this section, we recall some facts from algebraic integral geometry which will be needed later on. 

Let $V$ be a euclidean vector space of dimension $n$. Denote by  $SV=V \times S^{n-1}$ its sphere bundle and by $\pi=\pi_1:SV \to V, \pi_2:SV \to S^{n-1}$ the projections. We set $\Omega^{k,l}$ for the space of translation invariant differential forms of bidegree $(k,l)$ on $SV$. 

Using coordinates $(x_1,\ldots,x_n)$ on $V$, and $y_1,\ldots,y_n$ with $\sum_i y_i^2=1$ on $S^{n-1}$, the canonical $1$-form on $SV$ is 
\begin{displaymath}
\alpha=\sum_i y_i \ed x_i \in \Omega^{1,0}.
\end{displaymath}
The corresponding Reeb vector field is
\begin{displaymath}
T=\sum_i y_i \frac{\partial}{\partial x_i}.
\end{displaymath}

\begin{Definition}
A smooth translation invariant valuation of degree $0 \leq k \leq n-1$ is a map of the form 
\begin{displaymath}
\phi(K)=\int_{\nc(K)}\omega, \quad K \in \mathcal{K}(V),
\end{displaymath}
where $\omega \in \Omega^{k,n-k-1}$ and where $\nc(K)$ denotes the normal cycle of $K$. A smooth translation invariant valuation of degree $n$ is of the form $\phi(K)=c \vol(K)$ with $c \in \R$. 

A smooth area measure of degree $0 \leq k \leq n-1$ is of the form 
\begin{displaymath}
\Phi(K,A)=\int_{\nc(K) \cap \pi_2^{-1}A}\omega, \quad K \in \mathcal{K}(V), A \in \mathcal{B}(S^{n-1}).
\end{displaymath}
The spaces of smooth translation invariant valuations and area measures of degree $k$ are denoted by $\Val_k^{sm}, 0 \leq k \leq n$ and $\Area_k^{sm}, 0 \leq k \leq n-1$.
\end{Definition}

The form $\omega$ in the definition is not unique. To describe its kernel, we need the Rumin operator $\D$. Given $\omega \in \Omega^{k,n-k-1}$, there exists a unique form $\alpha \wedge \xi \in \Omega^{k,n-k-1}$ such that $\ed(\omega+\alpha \wedge \xi)$ is divisible by $\alpha$ and $\D\omega:=\ed(\omega+\alpha \wedge \xi)$ \cite{rumin94}. We refer to the next section for more information on this second order differential operator. 
 
\begin{Theorem}[\cite{bernig_broecker07}] \label{thm_kernel}
 \begin{enumerate}
  \item The valuation $\phi$ is trivial if and only if $\D\omega=0$ and $\pi_*\omega=0$.
  \item The area measure $\Phi$ is trivial if and only if $\omega$ is in the ideal generated by $\alpha$ and $\ed \alpha$. 
 \end{enumerate}
\end{Theorem}

There is an obvious map $\glob:\Area_k^{sm} \to \Val_k^{sm}$, given by $\glob \Phi(K):=\Phi(K,S^{n-1})$. This is a special case of a moment map. 

\begin{Definition}
The $r$-th moment map, $r \geq 0$, is the map
\begin{displaymath}
M^r:\Area_k^{sm} \to \Val_k^{sm} \otimes \Sym^rV
\end{displaymath}
defined by 
\begin{displaymath}
M^r\Phi(K):=\int_{S^{n-1}} \vec{y}^r \ed\Phi(K,y) \in \Sym^rV.
\end{displaymath}
\end{Definition}

The first variation map $\delta:\Val_{k+1}^{sm} \to \Area_k, 0 \leq k \leq n-1$ is defined by the equation 
\begin{displaymath}
 \left.\frac{d}{dt}\right|_{t=0} \mu(K+tL)=\int_{S^{n-1}} h_L \ed(\delta \phi(K)),
\end{displaymath}
where $K,L$ are compact convex sets and where $h_L$ is the support function of $L$. The existence of such a map was shown by Wannerer \cite[Proposition 2.2]{wannerer_unitary_module}. It is not difficult so show that $M^1 \circ \delta=0$. 

In terms of differential forms, the map $\delta_k:\Val^{sm}_{k+1} \to \Area_k, 0 \leq k \leq n-1$ can be described as follows. For $0 \leq k \leq n-2$, let $\omega \in \Omega^{k+1,n-k-2}$ represent $\phi \in \Val^{sm}_{k+1}$. Then $\delta_k \phi$ is represented by the form $i_T\D\omega \in \Omega^{k,n-k-1}$. If $k=n-1$, then $\delta_{n-1}(c\vol):=c S_{n-1}$, where $S_{n-1}$ is the surface area measure. 

In the proof of the main theorems, we will use the contraction map, as defined in \cite{wannerer_area_measures}. Let $a \in \Sym^r V, b \in \Sym^s V$ with $s \geq r$. Identify $a$ with an element $\alpha \in \Sym^r V^*$ by using the isomorphism $V \cong V^*$ induced by the scalar product. Then $\contr(a,b):=\alpha \lrcorner b \in \Sym^{s-r}V$, which is defined by 
\begin{displaymath}
 \langle \alpha \lrcorner b,\beta\rangle=\langle b,\alpha \cdot \beta\rangle, \quad \beta \in \Sym^{s-r} V^*.
\end{displaymath}

More explicitly, if $e_1,\ldots,e_n$ is an orthonormal basis of $V$ and 
$a=\sum a^{i_1 \ldots i_r} e_{i_1} \otimes \cdots \otimes e_{i_r}, b=\sum b^{j_1 \ldots j_s} e_{j_1} \otimes \cdots \otimes e_{j_s}$ with symmetric coefficients $a^{i_1 \ldots i_r}, b^{j_1 \ldots j_s}$, then 
\begin{displaymath}
 \contr(a,b)=\contr(b,a)=\sum a^{i_1\ldots i_r} b^{i_1 \ldots i_r j_1 \ldots j_{s-r}} e_{j_1} \otimes \ldots e_{j_{s-r}}.
\end{displaymath}

The contraction satisfies 
\begin{equation} \label{eq_contraction_formula1}
 \contr(a,\contr(b,c))=\contr(ab,c),
\end{equation}
where $a,b,c$ are symmetric tensors such that the rank of $c$ is at least the sum of the ranks of $a$ and $b$. In particular, 
\begin{equation} \label{eq_contraction_formula2}
 \contr(\vec{y}a,\vec{y}b)=\contr(a,b), \quad \vec{y} \in S^{n-1},
\end{equation}
which we will use several times. Another useful identity is 
\begin{equation} \label{eq_contraction_formula3}
\contr(\vec{y}^{r_1+r_2},a \cdot b)=\contr(\vec{y}^{r_1},a) \cdot \contr(\vec{y}^{r_2},b),
\end{equation}
where $a$ is of rank $r_1$ and $b$ is of rank $r_2$.

Given differential forms $\omega_1 \in \Omega^*(SV) \otimes \Sym^rV$, $\omega_2 \in \Omega^*(SV) \otimes \Sym^s V$, the contraction $\contr(\omega_1,\omega_2)$ is defined by taking the wedge product on the form part and the contraction on the tensor part.  

\section{The translation invariant Rumin-de Rham complex}
\label{sec_rumin_de_rham}

We first recall the construction of the Rumin complex from \cite{rumin94}. Let $(M,H)$ be a contact manifold of dimension $2n-1$. Locally, it can be described by a $1$-form $\alpha$ (called contact form) such that $H=\ker \alpha$ and $\alpha \wedge \ed\alpha^{n-1} \neq 0$. 

Let $\Omega^k:=\Omega^k(M)$ be the space of differential $k$-forms on $M$. A form $\omega \in \Omega^k(M)$ is called vertical, if $\omega|_H=0$, or equivalently if $\omega=\alpha \wedge \tau$ for some $\tau \in \Omega^{k-1}$. If a global contact form $\alpha$ is given, then the Reeb vector field $T$ is defined by the conditions $i_T\alpha=1, i_T\ed\alpha=0$. In this case, we call a form horizontal if $i_T\omega=0$. 

We refer to \cite{huybrechts05} for some basic notions in symplectic linear algebra. An $(n-1)$-form $\omega$ is called primitive if $\ed\alpha \wedge \omega$ is vertical. Given any $\omega \in \Omega^{n-1}$, there exists a form $\rho \in \Omega^{n-3}$ such that $\omega+\ed\alpha \wedge \rho$ is primitive.

Define two subspaces of $\Omega^k$ by 
\begin{align*}
\mathcal{I}^k & =\{ \omega \in \Omega^k: \omega= \alpha \wedge
\xi+ \ed\alpha \wedge \psi, \xi \in \Omega^{k-1}, \psi \in \Omega^{k-2}\}, \\ 
\mathcal{J}^k & =\{\omega \in \Omega^k: \alpha \wedge \omega=\ed\alpha
\wedge \omega=0\}.
\end{align*}

These spaces only depend on $(M,H)$ and not on the particular choice of a contact form $\alpha$. 

Since $\ed \mathcal{I}^k \subset \mathcal{I}^{k+1}$, there
exists an induced operator $\ed_H:\Omega^k/\mathcal{I}^k \rightarrow
\Omega^{k+1}/\mathcal{I}^{k+1}$. 

Similarly, $\ed \mathcal{J}^k \subset \mathcal{J}^{k+1}$ and the
restriction of $\ed$ to $\mathcal{J}^{k}$ yields an operator
$\ed_H:\mathcal{J}^k \rightarrow  \mathcal{J}^{k+1}$.  

In the middle dimension, there is a further operator, called {\it Rumin operator}, which is defined as follows. Let $\omega \in \Omega^{n-1}$. There exists $\xi \in
\Omega^{n-2}$ such that $\ed(\omega+\alpha \wedge \xi) \in
\mathcal{J}^n$, and this last form, which is unique, is denoted by $\D\omega$. It can be checked that $\D|_{\mathcal{I}^{n-1}}=0$, hence there
is an induced operator $\D:\Omega^{n-1}/\mathcal{I}^{n-1}
\rightarrow \mathcal{J}_{n}$. We will also need the operator $\Q:\Omega^{n-1} \to \Omega^{n-1}$ defined by $\Q\omega:=\omega+\alpha \wedge \xi$ (where $\xi$ is, as above, such that $\ed(\omega+\alpha \wedge \xi)$ is vertical). It satisfies $\Q^2=\Q, \ed \circ \Q=\D$.

The Rumin complex of the contact manifold $(M,H)$ is given by 

\begin{multline*}
0 \rightarrow C^\infty(M) \stackrel{\ed_H}{\rightarrow} \Omega^1/\mathcal{I}^1 \stackrel{\ed_H}{\rightarrow}
\ldots \stackrel{\ed_H}{\rightarrow} \Omega^{n-2}/\mathcal{I}^{n-2} \stackrel{\ed_H}{\rightarrow} \Omega^{n-1}/\mathcal{I}^{n-1}
\stackrel{\D}{\rightarrow} \mathcal{J}_{n} \stackrel{\ed_H}{\rightarrow} \\
\stackrel{\ed_H}{\rightarrow} \mathcal{J}_{n+1}
\stackrel{\ed_H}{\rightarrow} \ldots \stackrel{\ed_H}{\rightarrow}
\mathcal{J}_{2n-1} \rightarrow 0.
\end{multline*}

The cohomology of this complex is called {\it Rumin cohomology} and
denoted by $H_{\text{Rum}}^*(M,\R)$. By \cite{rumin94}, there exists a natural
isomorphism between Rumin cohomology and de Rham cohomology:  
\begin{equation} \label{rumin_isomorphism}
H_{\text{Rum}}^*(M,\R) \stackrel{\cong}{\longrightarrow} H^*_{\text{dR}}(M,\R).
\end{equation}

We will apply this theorem in the special case $M=SV$, where $V$ is a euclidean vector space of dimension $n$. 

Let $\Omega^{k,l}:=\Omega^{k,l}(SV)^{tr}$ denote the space of translation invariant forms on $SV$ of bidegree $(k,l)$. We use the convention that $\Omega^{k,l}=0$ if $k<0$ or $l<0$ or $k>n$ or $l>n-1$. Note that $\ed\Omega^{k,l} \subset \Omega^{k,l+1}$. Similarly, we define $\mathcal{I}^{k,l}, \mathcal{J}^{k,l}$ as the corresponding spaces of translation invariant forms of given bidegree.

The above Rumin-de Rham complex may be refined, for each $0 \leq k \leq n-1$, as  
\begin{multline*}
0 \rightarrow \Omega^{k,0}/\mathcal{I}^{k,0} \stackrel{\ed_H}{\rightarrow} \Omega^{k,1}/\mathcal{I}^{k,1} \stackrel{\ed_H}{\rightarrow}
\ldots \stackrel{\ed_H}{\rightarrow} \Omega^{k,n-k-1}/\mathcal{I}^{k,n-k-1} 
\stackrel{\D}{\rightarrow} \mathcal{J}^{k,n-k} \stackrel{\ed_H}{\rightarrow} \\
\stackrel{\ed_H}{\rightarrow} \mathcal{J}^{k,n-k+1}
\stackrel{\ed_H}{\rightarrow} \ldots \stackrel{\ed_H}{\rightarrow}
\mathcal{J}^{k,n-1} \rightarrow 0.
\end{multline*}

We call this complex the translation invariant Rumin-de Rham complex of degree $k$. 

\begin{Proposition} \label{prop_coh_translationinvariant}
 The cohomology of the translation invariant Rumin-de Rham complex vanishes, except at $\Omega^{k,0}/\mathcal{I}^{k,0}$ and at $\mathcal{J}^{k,n-1}$, where the cohomology is isomorphic to $\largewedge^kV^*$. 
\end{Proposition}

The proof is analogous to \cite{rumin94}, with \cite[Lemma 2.5]{bernig_qig} replacing the Poincar\'e lemma.

Let us describe the non-zero part of the cohomology more precisely. If $\omega \in \Omega^{k,0}$ with $\ed_H\omega=0$, then we may write $\ed \omega=\alpha \wedge \xi+\ed \alpha \wedge \psi$. We claim that $\omega':=\omega-\alpha \wedge \psi$ is $\ed$-closed. Indeed, $\ed\omega'=\alpha \wedge (\xi+\ed\psi)$. Taking $\ed$ once more yields $\ed \alpha \wedge (\xi+\ed \psi) \equiv 0 \mod \alpha$. Symplectic linear algebra tells us that $\xi+\ed \psi \equiv 0 \mod \alpha$, and hence $\ed \omega'=0$. It is easy to check that $\omega'$ only depends on the class of $\omega$ in $\Omega^{k,0}/\mathcal{I}^{k,0}$. 

Let us write 
\begin{displaymath}
\omega'=\sum_I f_I \kappa_I,
\end{displaymath}
where $\kappa_I$ ranges over a basis of $\largewedge^k V^*$ and where $f_I$ is a function on the sphere. Since $\ed\omega'=0$, the functions $f_I$ are in fact constants, and we map the Rumin cohomology class $[\omega]$ to $\sum_I f_I \kappa_I \in \largewedge^k V^*$.

Next, let $\omega \in \mathcal{J}^{k,n-1}$. We may write 
\begin{displaymath}
\omega=\sum_I \kappa_I \wedge \tau_I,
\end{displaymath}
where $\kappa_I$ ranges over a basis of $\largewedge^k V^*$ and $\tau_I \in \Omega^{n-1}(S^{n-1})$. Now $\tau_I$ is exact if and only if $\int_{S^{n-1}} \tau_I=0$, and the map which sends the cohomology class $[\omega]$ to $\pi_*\omega=\sum_I \kappa_I \int_{S^{n-1}} \tau_I \in \largewedge^kV^*$ establishes an isomorphism.

\begin{Proposition} \label{prop_rumin_multiple}
Let $V$ be a euclidean vector space of dimension $n$. Let $\omega \in \Omega^{k,n-k-1}$ be primitive (i.e. $\ed \alpha \wedge \omega=\alpha \wedge \tau$ for some $\tau$), and such that $\D\omega=\ed\omega$. Let $f \in C^\infty(S^{n-1})$. Define a translation invariant vector field $X_f$ on $SV$ by the conditions 
\begin{displaymath}
i_{X_f}\alpha=0 \quad i_{X_f}\ed\alpha=df.
\end{displaymath}
Then 
\begin{displaymath}
\D(f\omega)=\ed(f \omega+\alpha \wedge \xi)=f\D\omega-\alpha \wedge (i_{X_f}\tau+\ed \xi),
\end{displaymath}
where $\xi:=i_{X_f}\omega \in \Omega^{k,n-k-2}$. 
\end{Proposition}

\proof
The existence of $X_f$ follows by basic properties of contact manifolds: the first condition means that $X_f$ belongs to the contact plane. Since $\ed\alpha$ is non-degenerated on the contact plane, the second condition fixes $X_f$ uniquely. 

Let us write 
\begin{displaymath}
\ed\alpha \wedge \omega=\alpha \wedge \tau
\end{displaymath}
for some form $\tau \in \Omega^{k,n-k}$. This is possible by the assumption that $\omega$ is primitive. Plugging $X_f$ into this equation yields
\begin{displaymath}
\ed f \wedge \omega+\ed\alpha \wedge \xi=-\alpha \wedge i_{X_f}\tau.
\end{displaymath}
Now we compute
\begin{align*}
\ed(f\omega+\alpha \wedge \xi) & =\ed f \wedge \omega+f \ed\omega+\ed\alpha \wedge \xi-\alpha \wedge \ed\xi \\
& = -\alpha \wedge i_{X_f}\tau+f\ed\omega-\alpha \wedge \ed \xi. 
\end{align*}
This form is vertical, since $\ed\omega$ is vertical by assumption. 
\endproof


\section{The kernel of the centroid map}
\label{sec_centroid_map}

\proof[Proof of Theorem \ref{mainthm_moments}, part (1)]

It is well-known and easy to prove that the image of $\delta_k$ is in the kernel of $M^1$ for each $0 \leq k \leq n-1$, see \cite[Lemma 2.9]{wannerer_unitary_module}. We have to prove the reverse inclusion. 

\begin{itemize}
\item[Case $k=0$]
Let $\Phi \in \Area_0$ be represented by the horizontal form $\omega \in \Omega^{0,n-1}$ and belong to the kernel of $M^1$. Then $\vec{y} \omega$ defines the trivial tensor valuation, which implies by Theorem \ref{thm_kernel} that $\int_{S^{n-1}} \vec{y}\omega=0$. We thus find $\tau_i \in \Omega^{n-2}(S^{n-1})$ such that $\ed\tau_i=y_i \omega$. Then 
\begin{displaymath}
 \alpha \wedge \omega=\sum_i \ed x_i \wedge y_i\omega=\sum_i \ed x_i \wedge \ed\tau_i=-\ed\left(\sum_i \ed x_i \wedge \tau_i\right)=\D \phi,
\end{displaymath}
where $\phi:=-\sum_i \ed x_i \wedge \tau_i$. Contracting with $T$ yields $\omega=i_T\D\phi$, hence $\Phi$ is in the image of the first variation map $\delta_0:\Val_1 \to \Area_0$.

 \item[Case $2 \leq k \leq n-2$] We will construct a map $L^1:\Val_k^{sm} \otimes V \to \Area^{sm}_k/\mathrm{im} \delta_k$ such that $L^1 \circ M^1$ equals the projection from $\Area^{sm}_k$ to $\Area^{sm}_k/\mathrm{im} \delta_k$. The corresponding diagram is 
\begin{displaymath}
 \xymatrix{
 0 \ar[r] & \Val_{k+1} \ar[r]^{\delta_k} & \Area_k \ar[r]^{M^1} \ar[dr] & \Val_k \otimes V \ar[d]^{L^1} \\ & & & \Area_k/\mathrm{im}\delta_k
 }
\end{displaymath}

Obviously the existence of such a map implies the statement.

Let $\phi \in \Val_k \otimes V$ be represented by the form $\hat \omega \in \Omega^{k,n-k-1} \otimes V$. Set 
\begin{displaymath}
 \rho:=i_T\D\hat \omega \in \Omega^{k-1,n-k} \otimes V,
\end{displaymath}
which only depends on $\phi$ but not on the choice of $\hat \omega$. 

Define 
\begin{displaymath}
 \tau:=\frac{1}{k+1} \contr\left(\sum_j \ed x_j \cdot e_j,\rho\right) \in \Omega^{k,n-k}.
\end{displaymath}

Since $\ed(\alpha \wedge \rho)=\ed \D\hat \omega=0$ and $\ed\alpha \wedge \alpha \wedge \rho=\ed\alpha \wedge \D\hat \omega=0$, we obtain that $\ed(\alpha \wedge \tau)=0$ and  $\ed\alpha \wedge \alpha \wedge \tau=0$. Hence $\alpha \wedge \tau \in \mathcal{J}^{k+1,n-k}$ defines a Rumin cohomology class, which is trivial by Proposition \ref{prop_coh_translationinvariant} (here we use that $k \neq 1$). Hence $\alpha \wedge \tau=\ed\kappa$ for some $\kappa \in \mathcal{J}^{k+1,n-k-1}$.  

Put 
\begin{displaymath}
 \omega:=\Q i_T \kappa \in \Omega^{k,n-k-1}.
\end{displaymath}

We define $L^1\phi \in \Area_k/\mathrm{im}\delta_k$ to be the equivalence class of the area measure presented by $\omega$. 

Let us first check that $L^1\phi$ is well-defined. The choice of $\kappa$ is not unique: by Proposition \ref{prop_coh_translationinvariant}, $\kappa$ may be replaced by $\kappa'=\kappa+\D\eta$ for some $\eta \in \Omega^{k+1,n-k-2}$ (here we use that $k \neq n-1$). 

We have 
\begin{displaymath}
 0=\ed \D\eta=\ed(\alpha \wedge i_T\D\eta)=\ed\alpha \wedge i_T\D\eta-\alpha \wedge \ed i_T\D\eta. 
\end{displaymath}
Contracting with $T$ gives us $\ed i_T\D\eta=\alpha \wedge i_T\ed i_T\D\eta$, hence $\ed i_T\D\eta$ is vertical and $Q i_T\D\eta=i_T\D\eta$. The area measure defined by $\omega':=\Q i_T\kappa'$ differs from the area measure defined by $\omega$ by the area measure defined by $i_T\D\eta$, which is in the image of $\delta_k$. 

Let us next prove that $L^1 \circ M^1:\Area_k \to \Area_k/\mathrm{im}\delta_k$ is just the projection map.

Let $\omega \in \Omega^{k,n-k-1}$ represent an area measure $\Phi$. Adding suitable multiples of $\alpha$ and $\ed \alpha$, we may assume that $\ed \omega=\D\omega$ and $\ed \alpha \wedge \omega=\alpha \wedge \tau$ for some $\tau \in \Omega^{k,n-k}$ with $i_T\tau=0$. Setting $\kappa:=\alpha \wedge \omega$ we find that $\kappa \in \mathcal{J}^{k+1,n-k-1}, \ed \kappa=\alpha \wedge \tau$ and $\omega=\Q i_T\kappa$. Since $M^1\Phi$ is represented by $\hat \omega:=\vec{y}\omega$, it only remains to prove that 
\begin{equation} \label{eq_tau}
\tau \equiv \frac{1}{k+1} \contr\left(\sum_j \ed x_j \cdot e_j,i_T\D(\vec{y}\omega)\right) \mod \alpha.
\end{equation}

To compute $\D(\vec{y}\omega)$, we use Proposition \ref{prop_rumin_multiple} and set $\xi:=i_{X_{\vec{y}}}\omega$ where
\begin{displaymath}
X_{\vec y}=\sum_{i=1}^n X_{y_i} e_i=T \vec{y}-\sum_j \frac{\partial}{\partial x_j}e_j.
\end{displaymath}

Then
\begin{align*}
\ed\xi & = \ed(i_{X_{\vec{y}}}\omega)\\
& = \ed\left(\vec{y} i_T\omega-\sum_j i_{\frac{\partial}{\partial x_j}}\omega e_j \right)\\
& = \ed\vec{y}  \wedge i_T\omega+\vec{y} \ed i_T\omega-\sum_j \ed i_{\frac{\partial}{\partial x_j}}\omega e_j.
\end{align*}

Since $\omega$ is translation invariant and $\ed\omega$ is vertical, we obtain that
\begin{align}
\ed i_{\frac{\partial}{\partial x_j}}\omega & = - i_{\frac{\partial}{\partial x_j}}\ed\omega\\
& =  - i_{\frac{\partial}{\partial x_j}}(\alpha \wedge i_T\ed\omega)\\
& \equiv -y_j i_T\ed \omega \mod \alpha. \label{eq_translation_inv_omega}
\end{align}

Plugging this into the above equation yields 
\begin{displaymath}
\ed\xi \equiv  \ed\vec{y} \wedge i_T\omega+\vec{y} \ed i_T\omega+\vec{y} i_T\ed\omega \mod \alpha.
\end{displaymath}

Proposition \ref{prop_rumin_multiple} gives
\begin{align*}
\D(\vec{y}\omega) & =  \vec{y}\ed\omega-\alpha \wedge \left(i_{X_{\vec{y}}}\tau+\ed\xi\right) \\
& = \alpha \wedge \left(\sum_j i_{\frac{\partial}{\partial x_j}} \tau e_j-\ed\vec{y} \wedge i_T\omega-\vec{y} \ed i_T\omega \right).
\end{align*}

Putting 
\begin{displaymath}
 A_1 := \sum_j i_{\frac{\partial}{\partial x_j}} \tau e_j, \quad A_2  := \ed\vec{y}  \wedge i_T\omega, \quad A_3:=\vec{y} \ed i_T\omega,
\end{displaymath}
we thus have $i_T\D(\vec{y} \omega) \equiv A_1-A_2-A_3 \mod \alpha$.

Using \eqref{eq_contraction_formula1}, \eqref{eq_contraction_formula2} and \eqref{eq_contraction_formula3} we find
\begin{align*}
 \mathrm{contr}\left(\sum_j \ed x_j e_j,A_1\right) & =\sum_j \ed x_j \wedge i_\frac{\partial}{\partial x_j}\tau=k\tau \\
 \mathrm{contr}\left(\sum_j \ed x_j e_j,A_2\right) & = \sum_j \ed x_j \wedge \ed y_j \wedge i_T\omega=-\ed\alpha \wedge i_T\omega=-\tau\\
\mathrm{contr}\left(\sum_j \ed x_j e_j,A_3\right) & = \sum_j y_j \ed x_j \wedge \ed i_T\omega=\alpha \wedge \ed i_T\omega.
\end{align*}

It follows that 
\begin{displaymath}
 \mathrm{contr}\left(\sum_j \ed x_j e_j,i_T\D(\vec{y}\omega)\right) \equiv (k+1)\tau \mod \alpha,
\end{displaymath}
as claimed.
\item[Case $k=n-1$]
We adapt the proof of the case $2 \leq k \leq n-2$. The condition $k \neq n-1$ was only used for well-definedness of $L^1$. Using the same notation, let us do this part separately. Define $\tau \in \Omega^{n-1,1}$ as before. Let $\kappa,\kappa' \in \mathcal{J}^{n,0}$ be such that $\alpha \wedge \tau=\ed\kappa=\ed\kappa'$. We may write $\kappa=h \ed x_1 \wedge \ldots \wedge \ed x_n, \kappa'=h' \ed x_1 \wedge \ldots \wedge \ed x_n$ with $h,h'\in C^\infty(S^{n-1})$. From $\ed(\kappa-\kappa')=0$ we infer that $c:=h-h'$ is the constant function. Set (again as before) $\omega:=\Q i_T\kappa, \omega':=\Q i_T\kappa'$. 

We have
\begin{align*}
\ed i_T(\kappa-\kappa') & =\ed\left(c \sum_j (-1)^{j+1} \ed x_1 \wedge \ldots \wedge \ed x_{j-1} \wedge y_j \wedge \ed x_{j+1} \wedge \ldots \wedge \ed x_n\right)\\
& = c \sum_j \ed x_1 \wedge \ldots \wedge \ed x_{j-1} \wedge \ed y_j \wedge \ed x_{j+1} \wedge \ldots \wedge \ed x_n.
\end{align*}
It follows that 
\begin{align*}
 \alpha \wedge \ed i_T(\kappa-\kappa') & = c \sum_k y_k \ed x_k \wedge \sum_j \ed x_1 \wedge \ldots \wedge \ed x_{j-1} \wedge \ed y_j \wedge \ed x_{j+1} \wedge \ldots \wedge \ed x_n\\
 & = -c \sum_j y_j \ed y_j \wedge \ed x_1 \wedge \ldots \wedge \ed x_n\\
 & = 0.
\end{align*}
Hence $\omega-\omega'=\Q i_T(\kappa-\kappa')=i_T(\kappa-\kappa')=ci_T(\ed x_1 \wedge \ldots \wedge \ed x_n)$. The area measure defined by $\omega-\omega'$ therefore equals $cS_{n-1}= \delta_{n-1} c\vol$.  
The rest of the proof is again the same as in the previous case. 

\item[Case $k=1$]
Again, we adapt the proof of the case $2 \leq k \leq n-2$ and use the same notations. The condition $k \neq 1$ was used in order to write $\alpha \wedge \tau =\ed\kappa$ for some $\kappa$. In the case $k=1$, we have $\alpha \wedge \tau \in \mathcal{J}^{2,n-1}$ and trivially $\ed(\alpha \wedge \tau)=0$. However, since the cohomology does not vanish at this spot, we have to verify that $\pi_*(\alpha \wedge \tau)=\int_{S^{n-1}} \alpha \wedge \tau =0 \in \largewedge^2V^*$ (see the discussion after Proposition \ref{prop_coh_translationinvariant}). 

Write $\hat \omega=\sum_i \hat \omega_i e_i, \rho=\sum_i \rho_i e_i$ with $\omega_i \in \Omega^{1,n-2}, \rho_i \in \Omega^{0,n-1}=\Omega^{n-1}(S^{n-1})$. Then $\rho_i=i_T\D\hat \omega_i$. In other words: the area measure represented by $\rho_i$ equals the first variation map $\delta_0:\Val_1 \to \Area_0$ applied to the area measure represented by $\hat \omega_i$. Hence its first moment map vanishes, which by Theorem \ref{thm_kernel} implies that $\int_{S^{n-1}} \vec{y} \rho_i=0$. 

Next, we have $\tau=\frac12 \sum_i \ed x_i \wedge \rho_i$ and therefore $\alpha \wedge \tau=\frac12 \sum_{i,k} \ed x_k \wedge \ed x_i y_k\rho_i$. Since $ \int_{S^{n-1}} y_k\rho_i=0$ the claim follows. The rest of the proof is again as before. 
\end{itemize}
\endproof

Let us give another interpretation of the theorem and its proof. Let $\Phi \in \Area_k$ be such that $M^1\Phi=0$. The theorem tells us that $\Phi=\delta_k\mu$ for some valuation $\mu \in \Val_{k+1}$. Let us construct $\mu$ explicitly.

Let $\Phi$ be represented by a form $\omega \in \Omega^{k,n-k-1}(SV)$. By changing $\omega$ by multiples of $\alpha$ and $\ed\alpha$, we may assume as above that $\D\omega=\ed\omega$ and $\ed\alpha \wedge \omega=\alpha \wedge \tau$ for some $\tau \in \Omega^{k,n-k}$. Since $M^1\Phi=0$, we have $\D(\vec{y}\omega)=0$ and \eqref{eq_tau} implies that $\tau \equiv 0 \mod \alpha$. We infer that 
\begin{displaymath}
 \ed(\alpha \wedge \omega)=\ed\alpha \wedge \omega-\alpha \wedge \ed\omega=\alpha \wedge \tau-\alpha \wedge \D\omega=0.
\end{displaymath}
If $k \neq 0,n-1$ then \cite[Lemma 2.5]{bernig_qig} implies that $\alpha \wedge \omega=\ed\phi$ for some $\phi \in \Omega^{k+1,n-k-1}(SV)$. Then $\omega \equiv i_T\D\phi \mod \alpha$, i.e. $\Phi$ is the first variation measure of the valuation $\mu$ represented by $\phi$.

If $k=0$ we have $\omega \in \Omega^{0,n-1}$ and from $M^1\Phi=0$ and Theorem \ref{thm_kernel} we obtain that $\int_{S^{n-1}} \vec{y}\omega=0$. Similarly as in the proof above this implies that $\alpha \wedge \omega$ is exact. 

If $k=n-1$, then $\alpha \wedge \omega=h \ed x_1 \wedge \ldots \wedge \ed x_n$ for some function $h \in C^\infty(S^{n-1})$. Since $\ed(\alpha \wedge \omega)=0$, this function equals a constant $c$. It follows that $\omega \equiv c i_T(\ed x_1 \wedge \ldots \wedge \ed x_n) \mod \alpha$, hence $\Phi$ is a multiple of the surface area measure $S_{n-1}$, which equals the first variation measure of the volume.  

As a final remark concerning part (1) of the theorem, we give an interpretation in terms of curvature measures. To each $\omega \in \Omega^{k,n-k-1}, 0 \leq k \leq n-1$ one can associate a valuation $\Psi$ with values in the space of signed measures on $V$ by 
\begin{displaymath}
\Psi(K,\beta):=\int_{\nc(K) \cap \pi^{-1} \beta} \omega, \quad \beta \in \mathcal{B}(V). 
\end{displaymath}
Such valuations are called smooth translation invariant curvature measures of degree $k$, the corresponding space is denoted by $\Curv_k$ (compare \cite{bernig_fu_solanes}). 

If we let $\Phi$ denote the area measure induced by $\omega$, then the map $\Area_k \to \Curv_k, \Phi \mapsto \Psi$ is a well-defined bijection \cite{wannerer_unitary_module}. 

There is a first variation map $\delta_k:\Val^{sm}_{k+1} \to \Curv_k$, see \cite{bernig_fu_hig}. Our theorem may be restated as follows:

{\it A curvature measure $\Psi \in \Curv_k$ is in the image of $\delta_k$ if and only if for each compact convex body $K$ with smooth boundary and outer unit normal vector field $n$ we have
\begin{displaymath}
\int_{\partial K} \langle \xi,n\rangle d\Psi(K,\cdot)=0 \quad \forall \xi \in V.
\end{displaymath}}

This follows from Theorem \ref{mainthm_moments} and \cite[Proposition 4.18]{wannerer_unitary_module}.

\section{Injectivity of higher moment maps}
\label{sec_higher_moments}

\proof[Proof of Theorem \ref{mainthm_moments}, part (2)]
Let $1 \leq k \leq n-1$ and $r \geq 2$. We will construct a map $L^r:\Val_k \otimes \Sym^rV \to \Area_k$ such that $L^r \circ M^r=\mathrm{id}$. Clearly this implies the injectivity of $M^r$. 

We write $\widehat{\frac{\partial}{\partial y_j}}$ for the orthogonal projection of the vector field $\frac{\partial}{\partial y_j}$ on $V \times V$ onto $SV=V \times S^{n-1}$. Explicitly we have $\widehat{\frac{\partial}{\partial y_j}}=\frac{\partial}{\partial y_j}-y_j \sum_a y_a \frac{\partial}{\partial y_a}$. 

Let us define maps 
\begin{displaymath}
 \Omega^{k-1,n-k} \otimes \Sym^rV \stackrel{\So}{\to} \Omega^{k,n-k} \otimes V \stackrel{\Ro}{\to} \Omega^{k,n-k-1}
\end{displaymath}

 by
\begin{align*}
 \So \hat \rho & := \mathrm{contr}\left(\sum_j \ed x_j e_j \cdot \vec{y}^{r-2},\hat \rho\right)\\
 \Ro\omega & :=\sum_j\contr(e_j,i_{\widehat{\frac{\partial}{\partial y_j}}} \omega).
\end{align*}

Note that $\Ro$ applied to a multiple of $\alpha$ will again be a multiple of $\alpha$.

Let $\phi \in \Val_k \otimes \Sym^rV$ be represented by the form $\hat \omega \in \Omega^{k,n-k-1} \otimes \Sym^rV$. Define 
\begin{displaymath}
 \omega:=- \frac{2}{r(r-1)k(k+1)}  \Ro \circ \So i_T\D(\hat \omega) \in \Omega^{k,n-k-1}
\end{displaymath}
and let $L^r \phi \in \Area_k$ be the area measure presented by $\omega$. This map is obviously well-defined and we have to show that $L^r \circ M^r=\mathrm{id}$. 

Let $\Phi \in \Area_k$ be given by the form $\omega \in \Omega^{k,n-k-1}$. Then $M^r\Phi$ is given by the form $\hat \omega:=\vec{y}^r \omega$. 

Without loss of generality, we may assume that $\ed\omega=\D\omega$ and that $\ed\alpha  \wedge \omega=\alpha \wedge \tau$ for some $\tau \in \Omega^{k,n-k}$ with $i_T\tau=0$. We apply Proposition \ref{prop_rumin_multiple} to the $\Sym^r V$-valued function $f:=\vec{y}^r$. We compute  
\begin{displaymath}
X_{\vec{y}^r}=r\vec{y}^rT-r\sum \frac{\partial}{\partial x_j}e_j \vec{y}^{r-1}
\end{displaymath}
and hence 
\begin{align*}
\ed\xi & = \ed(i_{X_{\vec{y}^r}}\omega)\\
& = r \ed\left(\vec{y}^r i_T\omega-\sum_j i_{\frac{\partial}{\partial x_j}}\omega e_j \vec{y}^{r-1}\right)\\
& = r^2 \ed\vec{y} \cdot \vec{y}^{r-1} \wedge i_T\omega+r \vec{y}^r \ed i_T\omega-r \sum_j \ed i_{\frac{\partial}{\partial x_j}}\omega e_j \vec{y}^{r-1}\\
& \quad -(-1)^n r(r-1)\sum_j i_{\frac{\partial}{\partial x_j}}\omega \wedge \ed\vec{y} e_j \vec{y}^{r-2}.
\end{align*}

Plugging \eqref{eq_translation_inv_omega} into the above equation yields 
\begin{align*}
\ed\xi & \equiv  r^2 \ed\vec{y} \cdot \vec{y}^{r-1} \wedge i_T\omega+r \vec{y}^r \ed i_T\omega+r\vec{y}^{r} i_T\ed \omega\\
& \quad +(-1)^{n-1} r(r-1)\sum_j i_{\frac{\partial}{\partial x_j}}\omega \wedge \ed\vec{y} e_j \vec{y}^{r-2} \mod \alpha.
\end{align*}

Proposition \ref{prop_rumin_multiple} gives us
\begin{align*}
\D \hat \omega & = \vec{y}^r \ed\omega -\alpha \wedge \left(i_{X_{\vec{y}^r}}\tau+\ed\xi\right) \\
& = \vec{y}^r \ed \omega + r \alpha \wedge \bigg( \sum_j i_{\frac{\partial}{\partial x_j}} \tau e_j \vec{y}^{r-1}-r \ed\vec{y} \cdot \vec{y}^{r-1} \wedge i_T\omega-\vec{y}^r \ed i_T\omega\\
& \quad -\vec{y}^r i_T\ed\omega+(-1)^n (r-1)\sum_j i_{\frac{\partial}{\partial x_j}}\omega \wedge \ed\vec{y} e_j \vec{y}^{r-2} \bigg) \\
& = (1-r) \vec{y}^r \ed \omega+r \alpha \wedge \bigg(\sum_j i_{\frac{\partial}{\partial x_j}} \tau e_j \vec{y}^{r-1}-r \ed\vec{y} \cdot \vec{y}^{r-1} \wedge i_T\omega\\ 
& \quad -\vec{y}^r \ed i_T\omega+(-1)^n (r-1)\sum_j i_{\frac{\partial}{\partial x_j}}\omega \wedge \ed\vec{y} e_j \vec{y}^{r-2} \bigg).
\end{align*}

Putting 
\begin{align*}
A_0 & := \vec{y}^r i_T \ed \omega\\
A_1 & := \sum_j i_{\frac{\partial}{\partial x_j}} \tau e_j \vec{y}^{r-1}\\
A_2 & := \ed\vec{y} \cdot \vec{y}^{r-1} \wedge i_T\omega\\
A_3 & := \vec{y}^r \ed i_T\omega\\
A_4 & := \sum_j i_{\frac{\partial}{\partial x_j}}\omega \wedge \ed\vec{y} e_j \cdot \vec{y}^{r-2} 
\end{align*}
we thus have 
\begin{equation} \label{eq_main_equation}
i_T\D(\hat \omega)=(1-r)A_0+rA_1-r^2 A_2-rA_3+(-1)^n r(r-1)A_4 \equiv 0 \mod \alpha.
\end{equation}

Taking into account \eqref{eq_contraction_formula1}, \eqref{eq_contraction_formula2}, \eqref{eq_contraction_formula3} we compute
\begin{align*}
\So A_0 & = \mathrm{contr}\left(\sum_j \ed x_j e_j ,\vec{y}^2i_T \ed \omega\right) \equiv 0 \mod \alpha,\\ 
\So A_1 & =  \contr\left(\sum_j \ed x_j e_j,\sum_j i_{\frac{\partial}{\partial x_j}}\tau e_j \cdot \vec{y}\right) \equiv \frac12 k \tau \vec{y} \mod \alpha, \\
 \So A_2 & =\mathrm{contr}\left(\sum_j \ed x_j e_j,\ed\vec{y} \cdot \vec{y}\right) \wedge i_T\omega & \\
 & \equiv -\frac12 \ed\alpha \wedge i_T\omega \cdot \vec{y} \mod \alpha,\\
 & \equiv -\frac12 \tau \cdot \vec{y} \mod \alpha,\\
\So A_3 & = \mathrm{contr}\left(\sum_j \ed x_j e_j ,\vec{y}^2 \ed i_T\omega\right) \equiv 0 \mod \alpha,\\ 
 \So A_4 & = \mathrm{contr}\left(\sum_j \ed x_j e_j,\sum_j i_{\frac{\partial}{\partial x_j}}\omega \wedge \ed\vec{y} e_j\right)\\
 & = \frac12\left(\sum \ed x_j \wedge i_{\frac{\partial}{\partial x_j}} \omega \wedge \ed\vec{y}+\ed x_a \wedge i_{\frac{\partial}{\partial x_j}} \omega \wedge \ed y_a e_j\right)\\
 & =\frac12 k \omega \wedge \ed\vec{y}+\frac12 (-1)^{n-1} \ed\alpha \wedge \sum_j i_{\frac{\partial}{\partial x_j}}\omega e_j.
\end{align*}

On the other hand, we have 
\begin{displaymath}
 \sum_j i_{\frac{\partial}{\partial x_j}}(\ed\alpha \wedge \omega) \cdot e_j=(-1)^n \omega \wedge \ed\vec{y}+\ed\alpha \wedge \sum_j i_{\frac{\partial}{\partial x_j}}\omega e_j
\end{displaymath}
and 
\begin{displaymath}
 \sum_j i_{\frac{\partial}{\partial x_j}}(\alpha \wedge \tau) \cdot e_j=\tau \vec{y}-\alpha \wedge \sum_j i_{\frac{\partial}{\partial x_j}}\tau e_j.
\end{displaymath}
Since the left hand sides agree (recall that $\ed \alpha \wedge \omega=\alpha \wedge \tau$), we obtain that 
\begin{displaymath}
  (-1)^{n-1} \ed\alpha \wedge \sum_j i_{\frac{\partial}{\partial x_j}}\omega e_j \equiv  (-1)^{n-1} \tau \vec{y} + \omega \wedge \ed\vec{y} \mod \alpha.
\end{displaymath}

It follows that 
\begin{displaymath}
 \So i_TD\hat \omega \equiv \frac{r(k+1)}{2} \tau \vec{y}+(-1)^n \frac{r(r-1)(k+1)}{2} \omega \wedge d\vec{y} \mod \alpha.
\end{displaymath}

Since 
\begin{align*}
\Ro(\tau \vec{y}) & =\sum_j i_{\widehat{\frac{\partial}{\partial y_j}}} \tau y_j=0,\\
\Ro(\omega \wedge \ed\vec{y}) & = \sum_j \left(i_{\widehat{\frac{\partial}{\partial y_j}}}\omega \wedge \ed y_j +(-1)^{n-1} \omega \wedge i_{\widehat{\frac{\partial}{\partial y_j}}}\ed y_j\right) \\
& = (-1)^n \sum_j \ed y_j \wedge i_{\widehat{\frac{\partial}{\partial y_j}}}\omega + (-1)^{n-1} (n-1) \omega \\
& = (-1)^n(n-k-1)\omega+(-1)^{n-1}(n-1)\omega\\
& =(-1)^{n+1}k\omega,
\end{align*}
we find 
\begin{align*}
\Ro \circ \So i_T\D \hat \omega & \equiv - \frac{r(r-1)k(k+1)}{2} \omega \mod \alpha.
\end{align*}
The definition of $L^r$ thus implies that $L^r \circ M^r \Phi=\Phi$.
\endproof


\section{Dual area measures and local additive kinematic formulas}
\label{sec_dual_area_measures}

Let $V$ be an oriented euclidean vector space of dimension $n$. The space $\Omega^{n-1}$ of translation invariant $(n-1)$-forms on $SV$ is endowed with its usual Fr\'echet topology of uniform convergence on compact subsets of all partial derivatives. We have a surjection  $\Omega^{n-1} \to \Area$ and endow the latter space with the quotient topology. 

Let $\Area^*$ be the dual space to $\Area$. We call elements of $\Area^*$ {\it dual area measures}. The globalization map $\glob:\Area \to \Val$ induces an inclusion $\glob^*:\Val^* \to \Area^*$.

Since the space $\Area$ is identified with the quotient of the space of translation invariant differential $(n-1)$-forms on the sphere bundle by the ideal generated by $\alpha$ and $\ed\alpha$, the dual space $\Area^*$ consists of all translation invariant $(n-1)$-currents which are Legendrian, i.e. vanish on the ideal $\langle \alpha,\ed\alpha\rangle$. The subspace $\Val^*$ is simply the space of such currents which are in addition closed (i.e. {\it cycles} in the terminology of geometric measure theory).  

Let $*_1:\Omega^*(SV)^{tr} \to \Omega^*(SV)^{tr}$ be the linear operator from \cite{bernig_fu06}. It is defined by the condition 
\begin{displaymath}
 *_1(\rho_1 \wedge \rho_2)=(-1)^{\binom{n-k}{2}} *\rho_1 \wedge \rho_2,
\end{displaymath}
where $\rho_1 \in \largewedge^kV^*, \rho_2 \in \Omega^*(S^{n-1})$ and where $*:\largewedge^kV^* \to \largewedge^{n-k}V^*$ denotes the usual Hodge star operator.

\begin{Lemma}
The space 
\begin{displaymath}
 \mathcal{J}^{n,tr}:=\{\tau \in \Omega^n(SV)^{tr}: \alpha \wedge \tau=0, \ed\alpha \wedge \tau=0\} 
\end{displaymath}
is closed under the operation 
\begin{displaymath}
 \tau_1 * \tau_2:= *_1^{-1}(*_1\tau_1 \wedge *_1\tau_2).
\end{displaymath}
\end{Lemma}

\proof
It is easily checked that $\alpha \wedge \tau=0 \iff i_T*_1\tau=0$ and $\ed\alpha \wedge \tau=0 \iff \mathcal{L}_T*_1\tau=0$. If $*_1\tau=*_1\tau_1 \wedge *_1\tau_2$ with $\tau_1,\tau_2 \in \mathcal{J}^{n,tr}$, then $i_T*_1\tau=(i_T*_1\tau_1) \wedge *_1\tau_2 \pm *_1\tau_1 \wedge (i_T*_1\tau_2)=0$ and $\mathcal{L}_T*_1\tau=(\mathcal{L}_T*_1\tau_1) \wedge *_1\tau_2+ *_1\tau_1 \wedge (\mathcal{L}_T*_1\tau_2)=0$.
\endproof

Given a translation invariant $(2n-1)$-form $\rho$ on $SV$, $\pi_*\rho \in \Omega^n(V)$ is translation invariant, hence a multiple of the volume form. We will denote the factor by $\int \rho$. Note that $\int \rho$ is independent of the orientation of $V$ (since the orientation affects the fiber integration but also the volume form). 
 
\begin{Definition} \label{def_smooth_dual_area}
 A dual area measure $L \in \Area^*$ is called smooth if there exists  $\tau \in \mathcal{J}^{n,tr}$ such that 
 \begin{displaymath}
  \langle L,\Phi\rangle = \int \omega \wedge \tau,
 \end{displaymath}
whenever $\omega \in \Omega^{n-1}(SV)^{tr}$ represents $\Phi \in \Area$. The space of smooth dual area measures is denoted by $\Area^{*,sm}$.
\end{Definition}

We note that $\omega$ is unique up to multiples of $\alpha$ and $\ed\alpha$. Since $\tau$ vanishes on such multiples, $\omega \wedge \tau$ is uniquely defined. Note also that, by Poincar\'e duality, for $L \in \Area^{*,sm}$, the form $\tau$ is unique. Since a change of orientation of $V$ results in changing the sign of the form $\omega$ representing a fixed smooth area measure, the same holds true for $\tau$.

\begin{Definition} \label{def_convolution_dual}
 Let $L_1,L_2 \in \Area^{*,sm}$ be represented by forms $\tau_1,\tau_2 \in \mathcal{J}^{n,tr}$. Then we define $L_1 * L_2 \in \Area^{*,sm}$ as the smooth dual area measure represented by $\tau_1 * \tau_2=*_1^{-1} (*_1\tau_1 \wedge *_1\tau_2) \in \mathcal{J}^{n,tr}$.
\end{Definition}

Let us check that $L_1 * L_2$ is independent of the choice of an orientation. Reversing the orientation of $V$ results in changing signs in $\tau_1,\tau_2$. Since $*_1$ also depends on the orientation and appears three times in the definition of $\tau_1 * \tau_2$, this form changes its sign as required.  

\begin{Lemma}
Let $G$ be a transitive group. Then the transposed of the inclusion $\Area^G \subset \Area$ restricts to a surjective map 
\begin{displaymath}
q_G:\Area^{*,sm} \to \Area^{G*}.
\end{displaymath}
\end{Lemma}

\proof
The map 
\begin{displaymath}
(\Omega^{n-1}/\mathcal I^{n-1})^G \otimes \mathcal J_n^G \to \R, \quad \omega \otimes \tau \mapsto \int \omega \wedge \tau
\end{displaymath}
is a non-degenerate pairing by Poincar\'e duality. In other words, given $l \in \Area^{G*}$ we find $\tau \in \mathcal J_n^G$ with $l(\Phi)=\int \omega \wedge \tau$ whenever $\Phi \in \Area^G$ is represented by $\omega$. Then the smooth dual area measure $L$ represented by $\tau$ restricts to $l$, as claimed. 
\endproof 

\proof[Proof of Theorem \ref{main_thm_ftaig}]
Let $A:\Area^G \to \Area^G \otimes \Area^G$ be the local additive kinematic formula and $A^*:\Area^{G*} \otimes \Area^{G*} \to \Area^{G*}$ its adjoint. We want to show that the following diagram commutes.
\begin{displaymath}
 \xymatrix{
 \Area^{*,sm} \otimes \Area^{*,sm} \ar[d]^{q_G \otimes q_G} \ar[r]^-{*} & \Area^{*,sm} \ar[d]^{q_G}\\
 \Area^{G*} \otimes \Area^{G*} \ar[r]^-{A^*} & \Area^{G*}
 }
\end{displaymath}

We need some notation from \cite{bernig_hug, wannerer_area_measures}. We let $\Val^{sm,r}:=\Val^{sm} \otimes \Sym^rV$ and $\Val^{r,G}:=(\Val \otimes \Sym^rV)^G$. Elements of $\Val^{sm,r}$ are called smooth tensor valuations of rank $r$. By \cite{bernig_hug}, there are additive kinematic formulas 
\begin{displaymath}
a^{r_1,r_2}:\Val^{r_1+r_2,G} \to \Val^{r_1,G} \otimes \Val^{r_2,G} 
\end{displaymath}
such that 
\begin{displaymath}
a^{r_1,r_2}(\Phi)(K,L)=\int_G (\mathrm{id}^{\otimes r_1} \otimes g^{\otimes r_2}) \Phi(K+g^{-1}L)\ed g.
\end{displaymath}

There is a natural perfect bilinear pairing $\Val^{sm,r} \otimes \Val^{sm,r} \to \R$: contract the tensor part and take the lowest degree part of the convolution of the valuation parts (see \cite{bernig_hug} for details). We thus obtain for transitive $G$ a map $\widehat{\pd}^r:\Val^{r,G} \to (\Val^{r,G})^*$.

In \cite[Theorem 3.2]{bernig_hug} and \cite[Prop. 4.7]{wannerer_area_measures} it was shown that the following diagram commutes
\begin{displaymath}
 \xymatrix{
 \Area^G \ar[r]^-A \ar[d]^-{M^{r_1+r_2}}  & \Area^G \otimes \Area^G \ar[d]^{M^{r_1} \otimes M^{r_2}} \\
 \Val^{r_1+r_2,G} \ar[r]^-{a^{r_1,r_2}} \ar[d]^{\widehat{\pd}^{r_1+r_2}} & \Val^{r_1,G} \otimes \Val^{r_2,G} \ar[d]^{\widehat{\pd}^{r_1} \otimes \widehat{\pd}^{r_2}}\\
  (\Val^{r_1+r_2,G})^* \ar[r]^-{c_G^*} & (\Val^{r_1,G})^* \otimes (\Val^{r_2,G})^*
 }
\end{displaymath}

Let us describe the map $\widehat{\pd}^r \circ M^r:\Area^G \to \left(\Val^{r,G}\right)^*$. Let $\Phi \in \Area^G$ be given by the differential form $\omega \in \Omega^{n-1}$. Let $\mu \in \Val^{r,G}$ be given by the form $\omega' \in \Omega^{n-1} \otimes \Sym^rV$. Then, by definition of $M^r$ and by \cite[Prop. 4.2.]{wannerer_area_measures}
\begin{equation} \label{eq_image_mr}
\langle \widehat{\pd}^r \circ M^r \Phi,\mu\rangle=\int \omega \wedge \contr(\vec{y}^r,\D\omega').
\end{equation}

Dualizing the above diagram yields the commutative diagram
\begin{displaymath}
 \xymatrix{
 \Area^{G*}  &  \Area^{G*} \otimes \Area^{G*} \ar[l]_-{A^*}\\
 \Val^{r_1+r_2,G}  \ar[u]^{(\widehat{\pd}^{r_1+r_2} \circ M^{r_1+r_2})^*} & \Val^{r_1,G} \otimes \Val^{r_2,G} \ar[l]_-{c_G}  \ar[u]_{(\widehat{\pd}^{r_1} \circ M^{r_1})^* \otimes (\widehat{\pd}^{r_2} \circ M^{r_2})^*}
 }
\end{displaymath}

Let $\mu_i \in \Val^{r_i,G}, i=1,2$ be represented by forms $\omega_i' \in \Omega^{n-1} \otimes \Sym^{r_i}V$. By \eqref{eq_image_mr}, the dual area measure $L_i:=(\widehat{\pd}^{r_i} \circ M^{r_i})^* \mu_i$ is represented by the form $\tau_i:=\contr(\vec{y}^{r_i},\D\omega_i') \in \Omega^n$.  

By \cite[Eq. (43)]{wannerer_area_measures}, the convolution product $\mu_1 * \mu_2$ is represented by the form 
\begin{displaymath}
 \omega':=*_1^{-1}(*_1\omega_1' \wedge *_1\D\omega_2') \in \Omega^{n-1} \otimes \Sym^{r_1+r_2}V,
\end{displaymath}
where we wedge the form part and take symmetric product in the tensor part. Consequently, the dual area measure 
\begin{displaymath}
L:=A^*(L_1 \otimes L_2)=(\widehat{\pd}^{r_1+r_2} \circ M^{r_1+r_2})^* (\mu_1 * \mu_2)
\end{displaymath}
is represented by the form
\begin{displaymath}
 \tau:=\contr(\vec{y}^{r_1+r_2},\D\omega').
\end{displaymath}

We have
\begin{displaymath}
 \D\omega'=*_1^{-1}(*_1\D\omega_1' \wedge *_1\D\omega_2') \in \Omega^{n} \otimes \Sym^{r_1+r_2}V.
\end{displaymath}

Since $*_1$ acts componentwise, we have
\begin{displaymath}
 *_1\tau_i=\contr(\vec{y}^{r_i},*_1\D\omega_i'),
\end{displaymath}
and hence, using \eqref{eq_contraction_formula3}, 
\begin{align*}
 *_1 \tau_1 \wedge *_1\tau_2 & =\contr(\vec{y}^{r_1},*_1\D\omega_1') \wedge \contr(\vec{y}^{r_2},*_1\D\omega_2')\\
 & = \contr(\vec{y}^{r_1+r_2},*_1\D\omega_1' \wedge *_1\D\omega_2')\\ 
 & = *_1\tau.
\end{align*}

By the definition of the convolution product on dual area measures we obtain that $L=L_1 * L_2$. 

Taking into account that $\widehat{\pd}^r$ is an isomorphism for each $r$, this shows that whenever $L_1 \in \Area^{G*}$ is in the image of $(M^{r_1})^*$ and $L_2 \in \Area^{G*}$ is in the image of $(M^{r_2})^*$, then $A^*(L_1 \otimes L_2)=L_1 * L_2$. 

Let us take $r_1=r_2:=2$. By Theorem \ref{mainthm_moments}, $M^2:\Area_k^G \to \Val_k^{2,G}$ is injective for $k \geq 1$. Since $\Area_0^G$ is spanned by the volume measure on the sphere whose second moment is non-zero, this map is still injective for $k=0$. Hence $(M^2)^*:(\Val^{2,G})^* \to \Area^{G*}$ is surjective, which implies that $A^*(L_1 \otimes L_2)=L_1 * L_2$ for all $L_1,L_2 \in \Area^{G*}$. This finishes the proof. 
\endproof


\section{Local additive kinematic formulas in hermitian integral geometry}

\label{sec_hermitian_case}

Let us first introduce some notation taken from \cite{bernig_fu_hig}. Let $V:=\C^n$ be a hermitian vector space of dimension $n$. We use coordinates $z_j=x_j+iy_j,j=1\ldots,n$ on $V$ and associated coordinates $\zeta_j=\xi_j+i\eta_j$ with $\sum_{j=1}^n (\xi_j^2+\eta_j^2)=1$ on $S^{2n-1}$.

Define $1$-forms 
\begin{align*}
\alpha & := \sum_{j=1}^n (\xi_j \ed x_j + \eta_j \ed y_j) \in \Omega^{1,0},\\
\beta & := \sum_{j=1}^n (\xi_j \ed y_j-\eta_j \ed x_j) \in \Omega^{1,0},\\
\gamma & :=  \sum_{j=1}^n (\xi_j \ed\eta_j-\eta_j \ed\xi_j) \in \Omega^{0,1}.
\end{align*}
Let 
\begin{displaymath}
\theta_1:=\ed\beta, \quad \theta_0:=\frac12 \ed\gamma, \quad \theta_2:=\sum_j \ed x_j \wedge \ed y_j.
\end{displaymath}
With the constant
\begin{equation} \label{eq_defcnkq}
c_{n,k,q}:=\frac{1}{q!(n-k+q)!(k-2q)!\omega_{2n-k}}
\end{equation}
we put
\begin{align*}
\beta_{k,q} & := c_{n,k,q} \beta \wedge \theta_0^{n-k+q} \wedge \theta_1^{k-2q-1} \wedge \theta_2^q, \quad \max\{0,k-n\} \leq q < \frac{k}{2}<n,\\
\gamma_{k,q} & := \frac{c_{n,k,q}}{2} \gamma \wedge \theta_0^{n-k+q-1} \wedge \theta_1^{k-2q} \wedge \theta_2^q, \quad \max\{0,k-n+1\} \leq q \leq \frac{k}{2}<n.
\end{align*}
Together with the symplectic form $\theta_s:=-\ed\alpha$, these forms generate the algebra of all $\overline{\mathrm U(n)}$-invariant forms on $SV$.

We let $B_{k,q}, \Gamma_{k,q}$ be the area measures induced by these forms and $B:=\mathrm{span}\{B_{k,q}\} \subset \Area^{\mathrm U(n)}, \Gamma:=\mathrm{span}\{\Gamma_{k,q}\} \subset \Area^{\mathrm U(n)}$.

The next proposition generalizes \cite[Lemma 6.1]{wannerer_area_measures}. 
\begin{Proposition} \label{prop_grading} 
Let $B^*=\mathrm{span}\{B_{k,q}^*\}, \Gamma^*=\mathrm{span}\{\Gamma_{k,q}\}$. Then 
\begin{displaymath}
 B^* * B^*=\{0\}, B^* * \Gamma^* \subset B^*, \Gamma^* * \Gamma^* \subset \Gamma^*.
\end{displaymath}
\end{Proposition}

In terms of local additive kinematic formulas, this means that
\begin{align*}
A(\Gamma_{k,q}) & \subset \Gamma \otimes \Gamma,\\
A(B_{k,q}) & \subset B \otimes \Gamma + \Gamma \otimes B.
\end{align*}

\proof
Let us write 
\begin{align*}
\tilde \theta_0 & := \theta_0,\\
\tilde \theta_1 & := \theta_1-\alpha \wedge \gamma,\\
\tilde \theta_2 & := \theta_2-\alpha \wedge \beta,\\
\tilde \theta_s & := \theta_s-\beta \wedge \gamma. 
\end{align*}

Let $(x,v) \in SV$ be a fixed point. Then 
\begin{displaymath}
T_{(x,v)}SV=T_xV \oplus T_vS^{2n-1}=\R x \oplus \R ix \oplus W \oplus \R ix \oplus W,
\end{displaymath}
where $W \subset V$ is the hermitian orthogonal complement of $\R x$. The space of $\overline{\mathrm{U}(n)}$-invariant forms can be identified with the space of $\mathrm{U}(n-1)$-invariant elements in $\largewedge^*T^*_{(x,v)}SV$, where $\mathrm{U}(n-1)$ is the stabilizer of $v$. We may decompose
\begin{displaymath}
(\largewedge^*T^*_{(x,v)}SV)^{\mathrm U(n-1)}=\largewedge^*(\R x \oplus \R ix \oplus \R ix) \otimes \largewedge^*(W \oplus W)^{\mathrm U(n-1)}.
\end{displaymath}
The first factor is generated by $\alpha,\beta,\gamma$, while the second factor is generated by $\tilde \theta_0,\tilde \theta_1,\tilde \theta_2,\tilde \theta_s$.

Since multiples of $\alpha$ and of $\theta_s$ induce the trivial area measure, any form $\beta \wedge \Theta$, where $\Theta$ is a polynomial in $\tilde \theta_0,\tilde \theta_1,\tilde \theta_2,\tilde \theta_s$ of total degree $2n-2$ induces an area measure in $B$. Similarly, any form $\gamma \wedge \Theta$ induces an area measure in $\Gamma$.

Let $L \in \Area^{\mathrm U(n)*}$ be represented by a form $\tau \in \mathcal{J}^{2n,tr}$ as in Definition \ref{def_smooth_dual_area}. Since $\alpha \wedge \tau=0$, we may write $\tau=\alpha \wedge \tilde \tau$ for some $(2n-1)$-form $\tilde \tau$ with $i_T\tilde \tau=0$. Since $\beta,\gamma$ are $1$-forms, while the forms $\tilde \theta$ are all $2$-forms, $\tilde \tau=\beta \wedge \Theta_1+\gamma \wedge \Theta_2$ for some polynomials $\Theta_1,\Theta_2$ in the $\tilde \theta$. 

If $L \in B^*$, then $L$ annihilates every $\Gamma_{k,q}$. It follows that 
\begin{displaymath}
\int \gamma \wedge \Theta \wedge \alpha \wedge (\beta \wedge \Theta_1+\gamma \wedge \Theta_2)  =0
\end{displaymath}
for every choice of $\Theta$. Poincar\'e duality implies that $\Theta_1=0$, hence $\tau$ is divisible by $\gamma$. Similarly, if $L \in \Gamma^*$, then $\tau$ is divisible by $\beta$. The reverse implications hold trivially.

Take now $L_1,L_2 \in B^*$, represented by forms $\tau_1,\tau_2$ as above. Then $\tau_1,\tau_2$ are both divisible by $\gamma$. Since $*_1$ act only on the base part of a form, $*_1\tau_1,*_1\tau_2$ are also divisible by $\gamma$ and hence $\tau_1 * \tau_2=*_1^{-1}(*_1\tau_1 \wedge *_1\tau_2)=0$. This shows that $B^* * B^* = \{0\}$. 

Next, take $L_1 \in B^*, L_2 \in \Gamma^*$. Then $\tau_1$ is divisible by $\gamma$, and then also $\tau_1 * \tau_2$ is divisible by $\gamma$, which shows that $L_1 * L_2 \in B^*$.

Finally, let $L_1,L_2 \in \Gamma^*$. Then $\tau_1,\tau_2$ are divisible by $\alpha \wedge \beta$, hence $*_1\tau_1,*_1\tau_2$ are in the algebra generated by $\gamma$ and the $\tilde \theta$'s. Since the degrees of $*_1\tau_1, *_1 \tau_2$ are even, $\gamma$ will in fact not appear. The wedge product $*_1\tau_1 \wedge *_1\tau_2$ is also in the algebra generated by the $\tilde \theta$'s, and taking $*_1^{-1}$ gives us a form which is divisible by $\beta$, i.e. an element of $\Gamma^*$.  
\endproof

\begin{Definition}
 On the algebra $\Omega^*(SV)^{\overline{\mathrm U(n)}}$ we introduce an algebra isomorphism $J$ which acts on the basis elements by 
 \begin{displaymath}
J\alpha=\alpha,  J\beta=\gamma, J\gamma=-\beta, J\tilde \theta_0=\tilde \theta_0,J\tilde \theta_1=\tilde \theta_1, J\tilde \theta_2=\theta_2, J \tilde \theta_s=\tilde \theta_s.
 \end{displaymath}
\end{Definition}

Since $J$ fixes $\alpha$ and $\ed\alpha$, $J$ induces complex structures on $\Area^{\mathrm U(n)}$ and on $\Area^{\mathrm U(n)*}$. Explicitly, 
 \begin{align*}
  JB_{k,q} & =\frac{2c_{n,k,q}}{c_{n,k-1,q}} \Gamma_{k-1,q}, \quad \max\{0,k-n\} \leq q < \frac{k}{2}<n,\\
  J\Gamma_{k,q} & =- \frac{c_{n,k,q}}{2c_{n,k+1,q}} B_{k+1,q}, \quad \max\{0,k-n+1\} \leq q \leq \frac{k}{2}<n,\\
  J^*B_{k,q}^* & = - \frac{c_{n,k-1,q}}{2c_{n,k,q}} \Gamma^*_{k-1,q}, \quad \max\{0,k-n\} \leq q < \frac{k}{2}<n,\\
  J^*\Gamma_{k,q}^* & = \frac{2c_{n,k+1,q}}{c_{n,k,q}} B^*_{k+1,q}, \quad \max\{0,k-n+1\} \leq q \leq \frac{k}{2}<n.
 \end{align*}

\begin{Lemma} \label{lemma_star1j}
Let $\tau \in (\mathcal{J}^{k,2n-k})^{\overline{\mathrm U(n)}}$ be divisible by $\beta$. Then 
  \begin{displaymath}
   *_1 J \tau=\beta \wedge \gamma \wedge *_1\tau.
\end{displaymath}
\end{Lemma}

\proof
We may assume that $\tau$ is of the form $\tau=\beta \wedge \rho_1 \wedge \rho_2$,
 with $\rho_1 \in \largewedge^{k-1}V^*, \rho_2 \in \Omega^{2n-k}(S^{2n-1})$.
Then 
\begin{align*}
 \beta \wedge \gamma \wedge *_1 \tau & = (-1)^{\binom{2n-k}{2}} \beta \wedge \gamma \wedge *(\beta \wedge \rho_1) \wedge \rho_2\\
 & = (-1)^{\binom{2n-k}{2}+k} \beta \wedge *(\beta \wedge \rho_1) \wedge \gamma \wedge \rho_2. 
\end{align*}
On the other hand side,
\begin{align*}
 *_1 J\tau & = *_1 (\gamma \wedge \rho_1 \wedge \rho_2)\\
 & = (-1)^{k-1} *_1 (\rho_1 \wedge \gamma \wedge \rho_2)\\
 & = (-1)^{\binom{2n-(k-1)}{2}+k-1} *\rho_1 \wedge \gamma \wedge \rho_2.
\end{align*}
The statement now follows from $\beta \wedge *(\beta \wedge \rho_1)=(-1)^{k-1}*\rho_1$ and a careful checking of signs. 
\endproof

\begin{Lemma}
\begin{enumerate}
\item $J^*(L_1 * L_2)=J^*L_1 * L_2=L_1 *  J^*L_2$ if $L_1,L_2 \in \Gamma^*$.
\item $J^*L_1 * J^*L_2=- L_1 * L_2$ if $L_1 \in \Gamma^*, L_2 \in B^*$.
\end{enumerate}
\end{Lemma} 

\proof
We use the notation from the proof of Proposition \ref{prop_grading}. Let $L \in \Area^{\mathrm U(n)*}$ be represented by a form $\tau \in \mathcal{J}^{2n}$. A small computation reveals that $J^*L$ is represented by $-J\tau$.  

Let $L_1,L_2 \in \Gamma^*$ be represented by $\tau_1,\tau_2$. Then $\tau_1$ and $\tau_2$ are divisible by $\beta$. Set $\tau:=*_1^{-1}(*_1\tau_1 \wedge *_1\tau_2)$, this form represents $L_1 * L_2$. Then $\tau$ is divisible by $\beta$ and Lemma \ref{lemma_star1j} implies that 
\begin{align*}
 - J\tau & = - *_1^{-1}(\beta \wedge \gamma \wedge *_1\tau)\\
 & = - *_1^{-1}((\beta \wedge \gamma \wedge *_1\tau_1) \wedge *_1\tau_2)\\
 & = *_1^{-1}(*_1(-J\tau_1) \wedge *_1\tau_2).
\end{align*}
The form on the left hand side represents $J(L_1 * L_2)$, while the form on the right hand side represents $J^*L_1 * L_2$.  This proves $J^*(L_1 * L_2)=J^*L_1 *L_2$. Changing the roles of $L_1$ and $L_2$ and using that the convolution product is commutative yields $J^*L_1 * L_2=L_1 * J^*L_2$. 

If $L_1 \in \Gamma^*, L_2 \in B^*$, then $J^*L_2 \in \Gamma^*$ and hence 
\begin{displaymath}
 J^*L_1 * J^*L_2=L_1 * (J^*)^2L_2=-L_1 * L_2.
\end{displaymath}

\endproof

\proof[Proof of Theorem \ref{mainthm_kinform_hermitian}]
Let $\max\{0,m-n+1\} \leq r \leq \frac{m}{2} < n$. By Proposition \ref{prop_grading}, we have $A(\Gamma_{m,r}) \in \Sym^2 \Gamma$, i.e. we may write 
\begin{displaymath}
 A(\Gamma_{m,r})=\sum \tilde c_{k,q,l,p}^{m,r} \Gamma_{k,q} \otimes \Gamma_{l,p}
\end{displaymath}
for some constants with $\tilde c_{k,q,l,p}^{m,r}=\tilde c_{l,p,k,q}^{m,r}$. Globalizing yields 
\begin{displaymath}
 a(\mu_{m,r})=\sum \tilde c_{k,q,l,p}^{m,r} \mu_{k,q} \otimes \mu_{l,p},
\end{displaymath}
and comparison of coefficients with the global additive kinematic formula shows that $\tilde c_{k,q,l,p}^{m,r}=c_{k,q,l,p}^{m,r}$ for all $k,q,l,p$. As a side remark, it follows from this argument that $c_{k,q,l,p}^{m,r}=0$ if $k-q=n$ or $l-p=n$. This seems to be a non-trivial new fact about the global formulas.

Next, let $\max\{0,m-n\} \leq r < \frac{m}{2}<n$. We want to show that 
\begin{displaymath}
 A(B_{m,r})=- (J \otimes \mathrm{id} + \mathrm{id} \otimes J) \circ A(J B_{m,r}).
\end{displaymath}
Both sides of the equation belong to $(B \otimes \Gamma+\Gamma \otimes B) \cap \Sym^2 \Area^{\mathrm U(n)}$. Therefore it is enough to check that $B_{k,q}^* \otimes \Gamma_{l,p}^*$ yields the same value on both sides. 
For this, we compute
\begin{align*}
 \langle A(B_{m,r}),B_{k,q}^* \otimes \Gamma_{l,p}^*\rangle & = \langle B_{m,r},B_{k,q}^* * \Gamma_{l,p}^*\rangle\\
  & = - \langle B_{m,r},J^*(J^* B_{k,q}^* * \Gamma_{l,p}^*)\rangle \\
   & = - \langle JB_{m,r},J^* B_{k,q}^* * \Gamma_{l,p}^*\rangle
\end{align*}
and 
\begin{align*}
 \langle (J \otimes \mathrm{id} & + \mathrm{id} \otimes J)A(J B_{m,r}),B_{k,q}^* \otimes \Gamma_{l,p}^*\rangle \\
 & =\langle A(JB_{m,r}),J^*B_{k,q}^* \otimes \Gamma_{l,p}^*+B_{k,q}^* \otimes J^*\Gamma_{l,p}^*\rangle\\
 & = \langle JB_{m,r},J^*B_{k,q}^* * \Gamma_{l,p}^*+\underbrace{B_{k,q}^* * J^*\Gamma_{l,p}^*}_{=0}\rangle.
\end{align*}
\endproof

\def\cprime{$'$}

\end{document}